\documentclass[12pt]{article}

\usepackage[geometry]{ifsym}
\usepackage{verbatim}
\usepackage{fancyhdr}
\usepackage{graphicx}
\usepackage{mathrsfs}
\usepackage{amssymb}
\usepackage{cite}
\usepackage{epsfig}
\usepackage{pst-poly}     
\usepackage{pst-all}      

\newtheorem{pps}{Proposition}[section]
\newtheorem{cor}{Corollary}[section]
\newtheorem{lem}{Lemma}[section]
\newtheorem{thm}{Theorem}[section]

\newtheorem{exap}{Example}[section]

\newenvironment{pf}[1][Proof]{\noindent\textbf{#1.} }{\hfill\rule{1mm}{2mm}}

\makeatletter \@addtoreset{equation}{section} \makeatother

\begin{document}

\title{Maximum Independent Sets Partition of $(n,k)$-Star Graphs\thanks{The work was supported by NNSF
of China (Nos. 11401004, 11371028, 11471016).}}
\author{Fu-Tao Hu\footnote{\ Correspondence to: F.-T. Hu; e-mail: hufu@ahu.edu.cn} \\
{\small School of Mathematical Sciences, Anhui University, Hefei, 230601, P.R. China}
}
\date{}
\maketitle

\begin{quotation}

\textbf{Abstract}: The $(n,k)$-star graph is a very important computer modelling.
The independent number and chromatic number of a graph are two important parameters in graph theory.
However, we did not know the values of this two parameters of the $(n,k)$-star graph
since it was proposed. In~\cite{wcz13}, Wei et. al. declared that they determined the independent number
of the $(n,k)$-star graph, unfortunately their proof is wrong.
This paper generalize their result and present a maximum independent sets partition of $(n,k)$-star graph.
From that we can immediately deduce the exact value of the independent number and chromatic number of $(n,k)$-star graph.



\end{quotation}

\section{Introduction}

For graph-theoretical notation and terminology not defined here we
follow \cite{x03}. In particular, let $G=(V,E)$ be a simple undirected
graph without loops and multi-edges, where $V=V(G)$ is the
{\it vertex set}, and $E=E(G)$ is the {\it edge set}.
If $xy\in E(G)$, we call two vertices $x$ and $y$ are {\it
adjacent}. For a vertex $x$, all the vertices
adjacent to it are the {\it neighbors} of $x$.

A subset $S$ of $V$ is said to be an {\it independent set} if no two of vertices are
adjacent in $S$ of a graph $G$. The cardinality of a maximum independent set in a graph $G$
is called the {\it independent number} of $G$ and is denoted by $\alpha(G)$.
Let $C$ be a set of $k$ colours. A {\it $k$-vertex-colouring} (simply a {\it $k$-colouring}),
is a mapping $c:\,V\rightarrow C$ such that any two adjacent vertices are assigned
the different colours of graph $G$. A graph is {\it $k$-colourable} if it has a $k$-colouring.
The {\it chromatic number}, which is denoted by $\chi(G)$, is the minimum $k$ for which graph $G$ is $k$-colourable.

As we known, the interconnection networks take an important part in the parallel computing/communication systems.
An interconnection network can be modeled by a graph where the the processors are the vertices and the edges are
the communication links.

In 1989, Akers and Krishnamurthy~\cite{ak89} introduced the $n$-dimensional star graph $S_n$, which
has superior degree and diameter compared to the hypercube as well as it is highly hierarchical and symmetrical~\cite{dt94}.
However, the vertex cardinality of the $n$-dimensional star is $n!$. The gap between $n!$ and $(n+1)!$ is very large when $S_n$ is extended to $S_{n+1}$.
Chiang and Chen~\cite{cc95} in 1995 generalized the star graph $S_n$ to the $(n,k)$-star graph,
which preserves many good properties of the star graph and has smaller scale. Since he $(n,k)$-star graph was introduced,
it has received great attention in the literature~\cite{ck01,cd08,cg10,cl02,cl13,cc95,cc98,dc14,hh03,hl06,lx14,ld04,ld08,sq09,wx12,xs10,yl10,z12} since then on.

The independent number and chromatic number of a graph are two important parameters in graph theory.
In~\cite{wcz13}, Wei et. al. declared that they determined the independent number
of the $(n,k)$-star graph, their result is right, unfortunately their proof is wrong.
In Section 3, we will show a counterexample of their result.
In Section 4, we will present a maximum independent sets partition
and determine the exact value of the independent number and chromatic number of $(n,k)$-star graph.


\section{Preliminary results}

We use $[n]$ to denote the set $\{1,2,\ldots,n\}$, where $n$ is a positive integer.
A \emph{permutation} of $[n]$ is a sequence of $n$ distinct symbols of $u_i\in [n]$,
$u_1u_2\ldots u_n$. The $n$-dimensional \emph{star network}, denoted by
$S_n$, is a graph with the vertex set $$V(S_n)=\{u_1u_2\ldots u_n:\, u_i\in [ n], u_i\ne u_j~{\rm for}~i\ne j \}.$$
The edges are specified as follows:

$E(S_n)$: $u_1u_2\ldots u_n$ is adjacent to $v_1v_2\ldots v_n$ iff there exists $i$ with $2\le i\le n$ such that $v_j =u_j$ for
$j \notin \{1, i\}$, $v_1 =u_i$ and $v_i =u_1$;\\
The Star graphs are vertex-transitive $(n-1)$-regular of order $n!$.

Let $n$ and $k$ be two positive integers with $k\in [n-1]$, and let
$\Gamma_{n,k}$ be the set of all $k$-permutations on $[n]$, that is,
$\Gamma_{n,k}=\{p_1p_2\ldots p_k:\,p_i\in [ n]~{\rm and}~p_i\ne p_j~{\rm for}~i\ne j\}$.
In 1995, Chiang and Chen~\cite{cc95} generalized the Star graph to $(n,k)$-Star graph denoted by $S_{n,k}$
with vertex set $V(S_{n,k})=\Gamma_{n,k}$.
The adjacency is defined as follows: $p_1p_2\ldots p_i\ldots p_k$ is adjacent to\\
(1) $p_ip_2\ldots p_1\ldots p_k$ where $2\le i\le k$;\\
(2) $xp_2\ldots p_i\ldots p_k$ where $x\in [ n] \setminus \{p_i:\,i\in [ k]\}$.

By definition, $S_{n,k}$ is $(n-1)$-regular vertex-transitive with $n!/(n-k)!$ vertices.
Moreover, $S_{n,n-1}\cong S_n$ and $S_{n,1}\cong K_n$ where $K_n$ is the complete graph with order $n$.

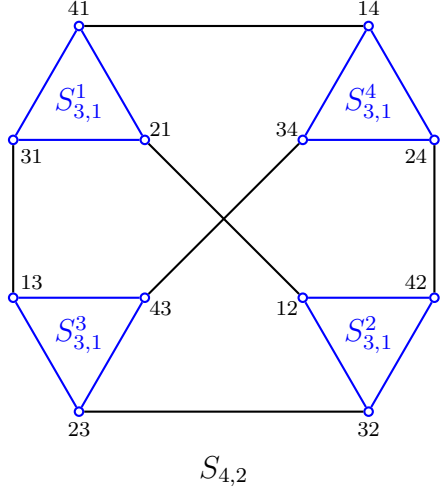
\begin{figure}[ht]
\psset{unit=0.7cm}
\begin{center}
\begin{pspicture}(-4.5,-5)(5,4)

\cnode[linecolor=blue](1.5,1.5){2pt}{41} \rput(1.2,1.7){\scriptsize 34}
\cnode[linecolor=blue](4,1.5){2pt}{42}  \rput(3.65,1.2){\scriptsize 24}
\cnode[linecolor=blue](2.75,3.665){2pt}{43}  \rput(2.75,4){\scriptsize 14}
\rput(2.75,2.2){\blue $S_{3,1}^4$}
\ncline[linecolor=blue]{41}{42}  \ncline[linecolor=blue]{42}{43}  \ncline[linecolor=blue]{43}{41}

\cnode[linecolor=blue](-1.5,1.5){2pt}{11} \rput(-1.2,1.7){\scriptsize 21}
\cnode[linecolor=blue](-4,1.5){2pt}{12}  \rput(-3.65,1.2){\scriptsize 31}
\cnode[linecolor=blue](-2.75,3.665){2pt}{13}  \rput(-2.75,4){\scriptsize 41}
\rput(-2.75,2.2){\blue $S_{3,1}^1$}
\ncline[linecolor=blue]{11}{12}  \ncline[linecolor=blue]{12}{13}  \ncline[linecolor=blue]{13}{11}

\cnode[linecolor=blue](-1.5,-1.5){2pt}{31} \rput(-1.2,-1.7){\scriptsize 43}
\cnode[linecolor=blue](-4,-1.5){2pt}{32}  \rput(-3.65,-1.2){\scriptsize 13}
\cnode[linecolor=blue](-2.75,-3.665){2pt}{33}  \rput(-2.75,-4){\scriptsize 23}
\rput(-2.75,-2.2){\blue $S_{3,1}^3$}
\ncline[linecolor=blue]{31}{32}  \ncline[linecolor=blue]{32}{33}  \ncline[linecolor=blue]{33}{31}

\cnode[linecolor=blue](1.5,-1.5){2pt}{21} \rput(1.2,-1.7){\scriptsize 12}
\cnode[linecolor=blue](4,-1.5){2pt}{22}  \rput(3.65,-1.2){\scriptsize 42}
\cnode[linecolor=blue](2.75,-3.665){2pt}{23}  \rput(2.75,-4){\scriptsize 32}
\rput(2.75,-2.2){\blue $S_{3,1}^2$}
\ncline[linecolor=blue]{21}{22}  \ncline[linecolor=blue]{22}{23}  \ncline[linecolor=blue]{23}{21}

\ncline{41}{31}  \ncline{42}{22}  \ncline{43}{13}  \ncline{11}{21}  \ncline{12}{32}  \ncline{23}{33}

\rput(0,-4.8){$S_{4,2}$}
\end{pspicture}
\caption{\label{f1}\footnotesize The $(4,2)$-star graph $S_{4,2}$.}
\end{center}
\end{figure}

Let $S_{n-1,k-1}^i$ denote a subgraph of $S_{n,k}$ induced by all the vertices with the same
last symbol $i$, for each $i\in [n]$. See Figure~\ref{f1} for instance.

\begin{lem}[Chiang and Chen~\cite{cc95}, 1995]
$S_{n,k}$ can be decomposed into $n$ subgraphs $S_{n,k}^i$, $i\in [n]$, and each subgraph
$S_{n-1,k-1}^i$ is isomorphic to $S_{n-1,k-1}$.
\end{lem}

\begin{lem}[Li and Xu~\cite{lx14}, 2014]\label{lem2.2}
For any $\alpha=p_2p_3\ldots p_k \in \Gamma_{n,k-1}~(k\ge 2)$, let $V_\alpha=\{p_1\alpha:\,p_1\in [n],~p_1\ne p_i,~2\le i\le k\}$.
Then the subgraph of $S_{n,k}$ induced by $V_\alpha$ is a complete graph of order $n-k+1$, denoted by $K_{n-k+1}^\alpha$.
\end{lem}

\section{Comments on Wei et. al.'s result}

In~\cite{wcz13}, Wei et. al. declared that they determined the independent number
of the $(n,k)$-star graph which was $\alpha(S_{n,k})=\frac{n!}{(n-k+1)!}$. 
We have checked their proof carefully. However, their proof is wrong.
They constructed an independent set $I_k$ of $S_{n,k}$ with cardinality $\frac{n!}{(n-k+1)!}$ step by step.
The following are their details.

Let $I_k^{(i,j)}$ be a vertex-set, which includes the vertices of $I_k$ if the vertices don't include
element $i$, and $I_k^{(i,j)}$ includes the vertices of $I_k$ if the vertices include element $i$ 
but swap $i$ with $j$.

\vskip6pt

Step 1: Clearly, $I_2=\{12,23,\ldots, (n - k +1) (n - k + 2), (n - k + 2)1\}$ is an independent set of $S_{n-k+2,2}$;

Step 2: In $S_{n-k+3,3}$, let $I_{3(n-k+3)}=\{\beta(n-k+3):\,\beta\in I_2\}$
and $I_{3x}=\{\beta x:\,x\in [n-k+2], \beta\in I_2^{(x,n-k+3)}\}$.
Now, we let  $I_3=\bigcup_{x\in[n-k+3]} I_{3x}$. Clearly, $|I_3|=(n-k+3)(n-k+2)$;

$\ldots \ldots \ldots$

{\bf Step k}: In $S_{n-k+3,3}$, let $I_{k(n-k+k)}=\{\beta(n-k+3):\,\beta\in I_{k-1}\}$
and $I_{kx}=\{\beta x:\,x\in [n-1], \beta\in I_{k-1}^{(x,n-k+k)}\}$.
Now, we let  $I_k=\bigcup_{x\in[n]} I_{kx}$. Clearly, $|I_k|=n(n-1)\ldots (n-k+3)(n-k+2)=\frac{n!}{(n-k+1)!}$.

The above $I_k$ is constructed in Wei's paper~\cite{wcz13}.
However, their construction is wrong and without detailed proof.
Next we use their method to construct the independent set $I_3$ for $S_{4,3}$.   

By Step 1, $I_2=\{12,23,31\}$ is a maximum independent set of $S_{3,2}$.
By Step 2, $I_{34}=\{124,234,314\}$, $I_{31}=\{421,231,341\}$, $I_{32}=\{142,432,312\}$, $I_{33}=\{123,243,413\}$.
Then $I_4=\bigcup_{x\in[4]} I_{3x}$. However, the vertices $124$ and $421$ are adjacent in $S_{4,3}$ (see Figure~\ref{f2}).
So their construction is wrong. 

In the next section, we will give a new method to construct a maximum independent sets partition of $S_{n,k}$ and show detailed proof for our result,
Which generalized Wei et. al.'s result.

\section{Maximum independent sets partition of $S_{n,k}$}

\begin{pps}\label{pps3.1}
The independent number of $S_{n,k}$ is $\alpha(S_{n,k})\le \frac{n!}{(n-k+1)!}$.
\end{pps}
\begin{pf}
This conclusion is true for $k=1$ since $S_{n,1}\cong K_n$. Next, assume $k\ge 2$.
Let $I$ be any maximum independent set of $S_{n,k}$.
For any $\alpha=p_2p_3\ldots p_k \in \Gamma_{n,k-1}~(k\ge 2)$, let $V_\alpha=\{p_1\alpha:\,p_1\in [n],~p_1\ne p_i,~2\le i\le k\}$.
Then the subgraph of $S_{n,k}$ induced by $V_\alpha$ is a complete graph of order $n-k+1$, denoted by $K_{n-k+1}^\alpha$ by Lemma~\ref{lem2.2}.
Thus, $I$ contains at most one vertex in $K_{n-k+1}^\alpha$. By definition,
there are exactly $\frac{n!}{(n-k+1)!}$ such $K_{n-k+1}^\alpha$.
Therefore, $\alpha(S_{n,k})=|I|\le \frac{n!}{(n-k+1)!}$.
\end{pf}

\begin{pps}\label{pps3.2}
Let $I_1^1=\{1\}, I_2^1=\{2\}, \ldots, I_n^1=\{n\}$. Then $\{I_1^1,I_2^1,\ldots,I_n^1\}$ is a maximum independent sets partition of $S_{n,1}$.
\end{pps}

\begin{pps}\label{pps3.3} Let
$$
 \begin{array}{ll}
&I_1^2=\{21,32,43,\ldots,n(n-1),1n\},\\
&I_2^2=\{31,42,53,\ldots,1(n-1),2n\},\\
&\ldots\\
&I_{n-1}^2=\{n1,12,23,\ldots,(n-2)(n-1),(n-1)n\}.
\end{array}
$$
Then $\{I_1^2,I_2^2,\ldots,I_{n-1}^2\}$ is a maximum independent sets partition of $S_{n,2}$.
\end{pps}

\vskip8pt
For each $j\in [n-k+1]$, we use $I_j^2$ of $S_{n-k+2,2}$
to generate a maximum independent set $I_j^3$ of $S_{n-k+3,3}$.
Step by step, we generate a maximum independent set $I_j^k$ of $S_{n,k}$
in the following. For each $i\in [k]\setminus \{1,2\}$, $\pi \in I_j^{i-1}$ and $x\in [n-k+i-1]$,
denote by $\pi(x,n-k+i)$ be such a permutation that replace $x$ by $n-k+i$ if
$x\in \pi$ ($x\in \pi$ means $x$ is equal to some symbol in $\pi$), otherwise $\pi(x,n-k+i)=\pi$.
Let $\pi=p_1p_2\ldots p_i$ be any vertex in $S_{n-k+i,i}$.
Denote by $\pi'=p_2p_1\ldots p_i$ be the vertex by exchanging the first two symbols in $\pi$.

\vskip6pt

{\bf Step 1}: By Proposition~\ref{pps3.3}, denote by $I_j^2=\{(j+1)1,(j+2)2, \ldots, (n-k+2)(n-k+2-j),1(n-k+3-j),
2(n-k+4-j),\ldots, j(n-k+2)\}$ be a
maximum independent set of $S_{n-k+2,2}$ for each $j\in [n-k+1]$.

{\bf Step 2}: Let $I_j^3(n-k+3)=\{\pi^{'}(n-k+3):\,\pi\in I_j^2\}$
and $I_j^3(x)=\{\pi(x,n-k+3)x:\,\pi\in I_j^2\}$ for each $x\in [n-k+2]$.
Let $I_j^3=\bigcup_{x\in[n-k+3]} I_j^3(x)$.

$\ldots \ldots \ldots$

{\bf Step i-1}: Let $I_j^i(n-k+i)=\{\pi^{'}(n-k+i):\,\pi\in I_j^{i-1}\}$
and $I_j^i(x)=\{\pi(x,n-k+i)x:\,\pi\in I_j^{i-1}\}$ for each $x\in [n-k+i-1]$.
Let $I_j^i=\bigcup_{x\in[n]} I_j^i(x)$.

$\ldots \ldots \ldots$

{\bf Step k-1}: Let $I_j^k(n)=\{\pi^{'}n:\,\pi\in I_j^{k-1}\}$
and $I_j^k(x)=\{\pi(x,n)x:\,\pi\in I_j^{k-1}\}$ for each $x\in [n-1]$.
Let $I_j^k=\bigcup_{x\in[n]} I_j^k(x)$.   {\hfill\rule{1mm}{2mm}}

\vskip6pt
Let $j\in [n-k+1]$ and $i\in [k]\setminus \{1,2\}$. By the above construction, it is easy to see that $I_j^i(x)$ is one corresponding to one with $I_j^{i-1}$ for each $x\in [n-k+i]$ and $I_j^i(x)\cap I_j^i(x')=\emptyset$ for every $x'\in [n-k+i]\setminus \{x\}$.
We can easily show the following conclusion by induction on $i$.

\begin{pps}\label{pps3.4}
For each $j\in [n-k+1]$ and $i\in [k]\setminus \{1\}$, $|I_j^i|=\frac{(n-k+i)!}{(n-k+1)!}$,
$I_j^i\cap I_{j'}^i=\emptyset$ for any $j'\in [n-k+1]\setminus \{j\}$.
Therefore $\{I_1^i, I_2^i,\ldots, I_{n-k+1}^i\}$ is a vertex sets partition of $S_{n-k+i,i}$.
\end{pps}

In the following, we show that $I_j^i$  is an independent set of $S_{n-k+i,i}$ for each $i\in \{3,4,\ldots,k\}$, $k\ge 3$ and $j\in [n-k+1]$.

\begin{lem}\label{lem3.1}
Let $i\in \{3,4,\ldots,k\}$, $k\ge 3$ and $j\in [n-k+1]$. If $I_j^i$ is an independent set of $S_{n-k+i,i}$,
then $I_j^{i-1}$ is an independent set of $S_{n-k+i-1,i-1}$.
\end{lem}
\begin{pf}
By Proposition~\ref{pps3.3}, $I_j^{2}$ is an independent set of $S_{n-k+i-1,i-1}$ for each $j\in [n-k+1]$.
Next assume $i\ge 4$. Suppose to the contrary that $I_j^{i-1}$ is not an independent set of $S_{n-k+i-1,i-1}$.
Firstly, assume $p_1p_2\ldots p_s\ldots p_{i-1}$ and $p_sp_2\ldots p_1\ldots p_{i-1}$
be two adjacent vertices in $I_j^{i-1}$ of $S_{n-k+i-1,i-1}$.
If $s=2$, then $p_1(n-k+i)\ldots p_{i-1}p_2$ and $(n-k+i)p_1\ldots p_{i-1}p_2$
are belong to $I_j^i(p_2)\subseteq I_j^i$ by the construction of $I_j^i$,
the one to one correspondence is
$$
 \begin{array}{ll}
&p_1p_2\ldots p_{i-1}\leftrightarrow p_1(n-k+i) \ldots p_{i-1}p_2,\\
&p_2p_1\ldots p_{i-1}\leftrightarrow (n-k+i)p_1 \ldots p_{i-1}p_2.
\end{array}
$$
However $p_1(n-k+i)\ldots p_{i-1}p_2$ and $(n-k+i)p_1\ldots \ldots p_{i-1}p_2$
are adjacent in $S_{n-k+i,i}$, a contradiction.
If $s>2$, then $p_1(n-k+i)\ldots p_s \ldots p_{i-1}p_2$ and $p_s(n-k+i)\ldots p_1 \ldots p_{i-1}p_2$
are belong to $I_j^i(p_2)\subseteq I_j^i$ by the construction of $I_j^i$,
the one to one correspondence is
$$
 \begin{array}{ll}
&p_1p_2\ldots p_s\ldots p_{i-1}\leftrightarrow p_1(n-k+i)\ldots p_s \ldots p_{i-1}p_2,\\
&p_sp_2\ldots p_1\ldots p_{i-1}\leftrightarrow p_s(n-k+i) \ldots p_1\ldots p_{i-1}p_2.
\end{array}
$$
However $p_1(n-k+i)\ldots p_s \ldots p_{i-1}p_2$ and $p_s(n-k+i)\ldots p_1 \ldots p_{i-1}p_2$
are adjacent in $S_{n-k+i,i}$, a contradiction.

Secondly, assume $p_1p_2\ldots p_{i-1}$ and $p_ip_2\ldots p_{i-1}$
be two adjacent vertices in $I_j^{i-1}$ of $S_{n-k+i-1,i-1}$.
By the construction of $I_j^i$, $p_1(n-k+i) \ldots p_{i-1}p_2$ and $p_i(n-k+i) \ldots p_{i-1}p_2$
are belong to $I_j^i(p_2)\subseteq I_j^i$, the one to one correspondence is
$$
 \begin{array}{ll}
&p_1p_2\ldots p_{i-1}\leftrightarrow p_1(n-k+i) \ldots p_{i-1}p_2,\\
&p_ip_2\ldots p_{i-1}\leftrightarrow p_i(n-k+i)\ldots p_{i-1}p_2.
\end{array}
$$

However $p_1(n-k+i)\ldots p_s \ldots p_{i-1}p_2$ and $p_i(n-k+i)\ldots p_s \ldots p_{i-1}p_2$
are adjacent in $S_{n-k+i,i}$, a contradiction.
\end{pf}

\begin{lem}\label{lem3.2}
Let $i\in \{3,\ldots,k\}$, $k\ge 3$ and $j\in [n-k+1]$. If $I_j^{i-1}$ is an independent set of $S_{n-k+i-1,i-1}$,
then the two vertices $p\pi_1$ and $p\pi_2$ can not both belong to $I_j^i$ of $S_{n-k+i,i}$
where $\pi_1$ and $\pi_2$ are adjacent in $S_{n-k+i-1,i-1}$.
\end{lem}
\begin{pf}
By Proposition~\ref{pps3.3},
$I_j^2=\{(j+1)1,(j+2)2, \ldots, (n-k+2)(n-k+2-j),1(n-k+3-j),
2(n-k+4-j),\ldots, j(n-k+2)\}$ is an independent set of $S_{n-k+2,2}$.

We consider the case for $i=3$.
Suppose to the contrary that there exist two vertices $p\pi_1$ and $p\pi_2$ in $I_j^3$
but $\pi_1$ and $\pi_2$ are adjacent in $S_{n-k+2,2}$.

Assume $\pi_1=p_1p_2$ and $\pi_2=p_2p_1$. Then $pp_1p_2\in I_j^3(p_2)$ and $pp_2p_1\in I_j^3(p_1)$.
If $p=n-k+3$, then $p_2p_1$ and
$p_1p_2$ are in $I_j^{2}$ by the construction of $I_j^3$,
the one to one correspondence is
$$
 \begin{array}{ll}
&(n-k+3)p_1p_2 \leftrightarrow p_2p_1,\\
&(n-k+3)p_2p_1\leftrightarrow p_1p_2.
\end{array}
$$
However $p_2p_1$ and $p_1p_2$ are adjacent in $S_{n-k+2,2}$, a contradiction with $I_j^2$ is an independent set of $S_{n-k+2,2}$.
If $p<n-k+3$, then $pp_1$ and
$pp_2$ are two vertices in $I_j^{2}$ by the construction of $I_j^3$,
the one to one correspondence is
$$
 \begin{array}{ll}
&pp_1p_2 \leftrightarrow pp_1,\\
&pp_2p_1\leftrightarrow pp_2.
\end{array}
$$
A contradiction with the construction of $I_j^2$.

Now, assume $\pi_1=p_1p_2$ and $\pi_2=p_3p_2$. Then $pp_1p_2\in I_j^3(p_2)$ and $pp_3p_2\in I_j^3(p_2)$.
If $p=n-k+3$, then $p_2p_1$ and
$p_2p_3$ are in $I_j^{2}$ by the construction of $I_j^3$,
the one to one correspondence is
$$
 \begin{array}{ll}
&(n-k+3)p_1p_2 \leftrightarrow p_2p_1,\\
&(n-k+3)p_3p_2\leftrightarrow p_2p_3.
\end{array}
$$
A contradiction with the construction of $I_j^2$.
If $p<n-k+3$, then $pp_1$ and
$pp_3$ are two vertices in $I_j^{2}$ by the construction of $I_j^3$,
the one to one correspondence is
$$
 \begin{array}{ll}
&pp_1p_2 \leftrightarrow pp_1,\\
&pp_2p_1\leftrightarrow pp_2.
\end{array}
$$
A contradiction with the construction of $I_j^2$.

Therefore, the conclusion is true for $i=3$.

We prove this Lemma by induction on $i$. Assume that the induction hypothesis is true for $i-1$ with $i\ge 4$.
We prove the case for $i\ge 4$.
Assume $I_j^{i-1}$ is an independent set of $S_{n-k+i-1,i-1}$. Then $I_j^{i-2}$ is an independent set of $S_{n-k+i-2,i-2}$ by Lemma~\ref{lem3.1}.
Suppose to the contrary that there exist two vertices $p\pi_1$ and $p\pi_2$ in $I_j^i$
but $\pi_1$ and $\pi_2$ are adjacent in $S_{n-k+i-1,i-1}$.

Firstly, assume $\pi_1=p_2p_3\ldots p_s \ldots p_{i-1}p_i$ and $\pi_2=p_sp_3\ldots p_2 \ldots p_{i-1}p_i$.
Suppose $s=i$. If $p=n-k+i$, then $p_ip_2p_3\ldots p_{i-1}$ and
$p_2p_ip_3\ldots p_{i-1}$ are in $I_j^{i-1}$ by the construction of $I_j^i$,
the one to one correspondence is
$$
 \begin{array}{ll}
&(n-k+i)p_2p_3 \ldots p_{i-1}p_i\leftrightarrow p_ip_2p_3 \ldots p_{i-1},\\
&(n-k+i)p_ip_3\ldots p_{i-1}p_2\leftrightarrow p_2p_ip_3 \ldots p_{i-1}.
\end{array}
$$
However $p_ip_2p_3\ldots p_{i-1}$ and $p_2p_ip_3\ldots p_{i-1}$ are adjacent in $S_{n-k+i-1,i-1}$, a contradiction.

Now assume $p<n-k+i$, then $pp_2p_3\ldots p_{i-1}$ and
$pp_ip_3\ldots p_{i-1}$ are two vertices in $I_j^{i-1}$ by the construction of $I_j^i$,
the one to one correspondence is
$$
 \begin{array}{ll}
&pp_2p_3 \ldots p_{i-1}p_i\leftrightarrow pp_2p_3 \ldots p_{i-1},\\
&pp_ip_3\ldots p_{i-1}p_2\leftrightarrow pp_ip_3 \ldots p_{i-1}.
\end{array}
$$
However $p_2p_3\ldots p_{i-1}$ and $p_ip_3\ldots p_{i-1}$ are adjacent in $S_{n-k+i-2,i-2}$,
a contradiction with the induction hypothesis.

Next, suppose $s\ne i$. If $p=n-k+i$, then $p_ip_2p_3\ldots p_s \ldots p_{i-1}$ and
$p_ip_sp_3\ldots p_2 \ldots p_{i-1}$ are two vertices in $I_j^{i-1}$ by the construction of $I_j^i$.
Now assume $p<n-k+i$. Then $pp_2p_3\ldots p_s \ldots p_{i-1}$ and
$pp_sp_3\ldots p_2 \ldots p_{i-1}$ are two vertices in $I_j^{i-1}$ by the construction of $I_j^i$.
However, $p_2p_3\ldots p_s \ldots p_{i-1}$ and
$p_sp_3\ldots p_2 \ldots p_{i-1}$ are two adjacent vertices in $S_{n-k+i-2,i-2}$, a contradiction with the induction hypothesis.

Secondly, assume $\pi_1=p_2p_3\ldots p_{i-1}p_i$ and $\pi_2=p_{i+1}p_3\ldots p_{i-1} p_i$.
If $p=n-k+i$, then $p_ip_2p_3\ldots p_{i-1}$ and
$p_ip_{i+1}p_3\ldots p_{i-1}$ are two vertices in $I_j^{i-1}$ by the construction of $I_j^i$.
Now suppose $p< n-k+i$. Then $pp_2p_3\ldots p_{i-1}$ and
$pp_{i+1}p_3\ldots p_{i-1}$ are two vertices in $I_j^{i-1}$ by the construction of $I_j^i$,
However $p_2p_3\ldots p_{i-1}$ and
$p_{i+1}p_3\ldots p_{i-1}$ are two adjacent vertices in $S_{n-k+i-2,i-2}$, a contradiction with the induction hypothesis.

By the principle of induction, this Lemma completes.
\end{pf}

\begin{lem}\label{lem3.3}
Let $i\in \{3,\ldots,k\}$, $k\ge 3$ and $j\in [n-k+1]$.
If $I_j^{i-1}$ is an independent set of $S_{n-k+i,i}$, then $I_j^i(x)$ is an independent set
of $S_{n-k+i,i}$ for each $x\in [n-k+i]$.
\end{lem}
\begin{pf}
Assume $i=3$. Suppose to the contrary that $I_j^3(x)$ is not an independent set of $S_{n-k+3,3}$.
Assume $p_1p_2x$ and $p_2p_1x$ be two adjacent vertices in $I_j^3(x)$.
By the construction of $I_j^3(x)$, $p_1p_2,p_2p_1\in I_j^2$ if $n-k+3\notin \{p_1,p_2\}$
and $xp_2, p_2x \in I_j^2$ if $p_1=n-k+3$ (the case for $p_2=n-k+3$ is similar).
A contradiction with the construction of $I_j^2$ in Proposition~\ref{pps3.3}.
Assume $p_1p_2x$ and $p_3p_2x$ be two adjacent vertices in $I_j^3(x)$.
By the construction of $I_j^3(x)$, we have
   $$p_2p_1, p_2p_3 \in I_j^2, ~~~{\rm if}~x=n-k+3,$$
   $$p_1x, p_3x \in I_j^2, ~~~{\rm if}~p_2=n-k+3,$$
   $$xp_2, p_3p_2 \in I_j^2, ~~~{\rm if}~p_1=n-k+3,$$
   $$p_1p_2, xp_2 \in I_j^2, ~~~{\rm if}~p_3=n-k+3,$$
   $$p_1p_2, p_3p_2 \in I_j^2, ~~~{\rm otherwise}.$$
Any case above makes a contradiction with the construction of $I_j^2$ in Proposition~\ref{pps3.3}.

We proceed by induction on $i\ge 3$.
Assume that the induction hypothesis is true for $i-1$ with $i\ge 4$.
Next we prove $I_j^i(x)$ is an independent set
of $S_{n-k+i,i}$ for each $x\in [n-k+i]$. Suppose to the
contrary that $I_j^i(x)$ is not an independent set
of $S_{n-k+i,i}$.

Firstly, assume $p_1p_2\ldots p_s\ldots p_{i-1}x$ and $p_sp_2\ldots p_1\ldots p_{i-1}x$
be two adjacent vertices in $I_j^i(x)$. On one hand, suppose $x=n-k+i$.
By the construction of $I_j^i(n-k+i)$, we have
$$p_2p_1\ldots p_s\ldots p_{i-1}, p_2p_s\ldots p_1\ldots p_{i-1}\in I_j^{i-1},~~{\rm if}~s\ge 3,$$
$$p_2p_1\ldots p_{i-1},p_1p_2\ldots p_{i-1}\in I_j^{i-1},~~{\rm if}~s=2.$$
The first case makes a contradiction with the result in Lemma~\ref{lem3.2}, since
$p_1p_3\ldots p_s\ldots p_{i-1}$ and $p_sp_3\ldots p_1\ldots p_{i-1}$ are adjacent in $S_{n-k+i-2,i-2}$.
The second case makes a contradiction with $I_j^{i-1}$ is an independent set.
On the other hand, suppose $x<n-k+i$. By the construction of $I_j^i(x)$,
we have $p_1p_2\ldots p_s\ldots p_{i-1}, p_sp_2\ldots p_1\ldots p_{i-1}\in I_j^{i-1}$ if $p_t\ne n-k+i$ for each $t=1,2,\ldots,i-1$
(if $p_t=n-k+i$ for some $t\in \{1,2,\ldots,i\}$, we just replace $p_t$ by $x$). However, $p_1p_2\ldots p_s\ldots p_{i-1}$ and $p_sp_2\ldots p_1\ldots p_{i-1}$ are adjacent in $S_{n-k+i-1,i-1}$, a contradiction with $I_j^{i-1}$ is an independent set.

Secondly, assume $p_1p_2\ldots p_{i-1}x$ and $p_ip_2\ldots p_{i-1}x$ be two adjacent vertices in $I_j^i(x)$.
If $x=n-k+i$, then $p_2p_1\ldots  p_{i-1}$ and $p_2p_i\ldots p_{i-1}$ are two vertices in $I_j^{i-1}$ by the construction of $I_j^i(n-k+i)$.
However $p_1\ldots p_{i-1}$ and $p_i\ldots p_{i-1}$ are two adjacent vertices in $I_j^{i-2}$, a contradiction with the result in Lemma~\ref{lem3.2}.
Assume $x<n-k+i$. By the construction of $I_j^i(x)$,
we have $p_1p_2\ldots p_{i-1}, p_ip_2\ldots p_{i-1}\in I_j^{i-1}$ if $p_t\ne n-k+i$ for each $t=1,2,\ldots,i-1,i$
(if $p_t=n-k+i$ for some $t\in \{1,2,\ldots,i\}$, we just replace $p_t$ by $x$). However, $p_1p_2\ldots p_{i-1}$ and $p_ip_2\ldots p_{i-1}$ are adjacent in $S_{n-k+i-1,i-1}$, a contradiction with $I_j^{i-1}$ is an independent set.

By the principle of induction, this Lemma completes.
\end{pf}

\begin{thm}\label{thm3.1}
The vertex set $I_j^i$ is an independent set of $S_{n-k+i,i}$ for each $j\in [n-k+1]$ and $i\in \{2,3,\ldots,k\}$.
\end{thm}
\begin{pf}
We proceed by induction on $i\ge 2$. By Proposition~\ref{pps3.3},
$I_j^2$ is an independent set of $S_{n-k+2,2}$.
Assume that the induction hypothesis is true for $i-1$ with $i\ge
3$.

Assume $I_j^{i-1}$ is an independent set of $S_{n-k+i-1,i-1}$.
By Lemma~\ref{lem3.3}, $I_j^i(x)$ is an independent set of  $S_{n-k+i,i}$ for each $x\in [n-k+i]$.
Suppose to the contrary that $I_j^i$ is not independent.
Assume $p_1p_2\ldots p_{i-1}p_i$ and $p_ip_2\ldots p_{i-1}p_1$ be two
adjacent vertices in $I_j^i$ (This is the only possible case since $I_j^i(x)$ is an independent set of  $S_{n-k+i,i}$ for each $x\in [n-k+i]$).
If $p_1=n-k+i$, then $p_ip_2\ldots p_{i-1}$ and $p_2p_i\ldots p_{i-1}$
should be in $I_j^{i-1}$ by the construction of $I_j^i$,
but they are two adjacent vertices in $I_j^{i-1}$, a contradiction (the case for $p_i=n-k+i$ is similar).
If $p_s=n-k+i$ for some $s$ with $2\le s\le i-1$,
then $p_1p_2\ldots p_i \ldots p_{i-1}$ (replace $n-k+i$ by $p_i$) and $p_ip_2\ldots p_1 \ldots p_{i-1}$ (replace $n-k+i$ by $p_1$)
should be in $I_j^{i-1}$ by the construction of $I_j^i$,
but they are two adjacent vertices in $I_j^{i-1}$, a contradiction. Now assume $p_s\ne n-k+i$ for
each $s\in [i]$. Then $p_1p_2 \ldots p_{i-1}$ and $p_ip_2 \ldots p_{i-1}$
should be in $I_j^{i-1}$ by the construction of $I_j^i$,
but they are two adjacent vertices in $I_j^{i-1}$, a contradiction. This Theorem completes.
\end{pf}

\begin{thm}
The constructed $I_j^k$ is a maximum independent set
of $S_{n,k}$ for each $j\in [n-k+1]$ and $\alpha(S_{n,k})=\frac{n!}{(n-k+1)!}$.
Moreover, $\{I_1^k, I_2^k,\ldots, I_{n-k+1}^k\}$ is a maximum independent sets
partition of $S_{n,k}$.
\end{thm}
\begin{pf}
By Proposition~\ref{pps3.1}, $\alpha(S_{n,k})\le \frac{n!}{(n-k+1)!}$.
By Proposition~\ref{pps3.4}, $I_j^k=\frac{n!}{(n-k+1)!}$.
By Theorem~\ref{thm3.1}, $I_j^k$ is an independent set of $S_{n,k}$.
Therefore, $I_j^k$ for each $j\in [n-k+1]$ is a maximum independent set
of $S_{n,k}$ and $\alpha(S_{n,k})=\frac{n!}{(n-k+1)!}$. By Proposition~\ref{pps3.4},
$\{I_1^k, I_2^k,\ldots, I_{n-k+1}^k\}$ is a vertex sets
partition of $S_{n,k}$, so $\{I_1^k, I_2^k,\ldots, I_{n-k+1}^k\}$ is a maximum independent sets
partition of $S_{n,k}$.
\end{pf}

 \vskip8pt
Since $\{I_1^k, I_2^k,\ldots, I_{n-k+1}^k\}$ is a maximum independent sets
partition of $S_{n,k}$, we immediately obtain the chromatic number of $S_{n,k}$.

\begin{cor}
The chromatic number of $S_{n,k}$ is $\chi(S_{n,k})=n-k+1$.
\end{cor}

Next, we show the maximum independent sets partition of $S_{4,3}$ by our construction.

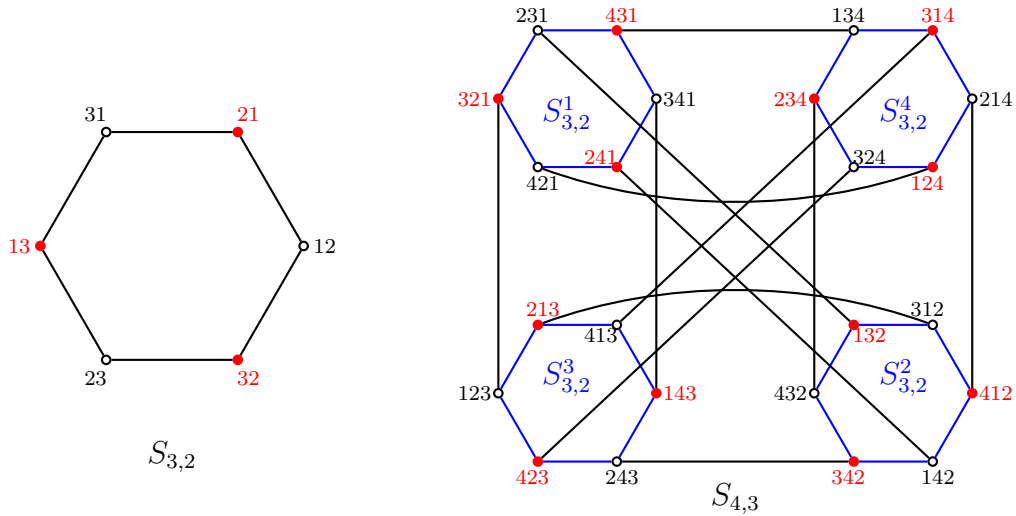
\begin{figure}[ht]
\psset{unit=0.7cm}
\begin{center}
\begin{pspicture}(-4.5,-5)(5,4.4)

\cnode(2.5;0){2pt}{1}\rput(2.9;0){\scriptsize 12}
\cnode*[linecolor=red](2.5;60){2pt}{2}\rput(2.9;60){\scriptsize\red 21}
\cnode(2.5;120){2pt}{3}\rput(2.9;120){\scriptsize 31}
\cnode*[linecolor=red](2.5;180){2pt}{4}\rput(2.9;180){\scriptsize\red 13}
\cnode(2.5;240){2pt}{5}\rput(2.9;240){\scriptsize 23}
\cnode*[linecolor=red](2.5;300){2pt}{6}\rput(2.9;300){\scriptsize\red 32}

\ncline{1}{2}  \ncline{2}{3}  \ncline{3}{4}  \ncline{4}{5}  \ncline{5}{6}  \ncline{6}{1}

\rput(0,-4){$S_{3,2}$}
\end{pspicture}
\begin{pspicture}(-5.5,-5)(5,4.4)

\cnode(2.25,1.5){2pt}{41} \rput(2.55,1.7){\scriptsize 324}
\cnode*[linecolor=red](3.75,1.5){2pt}{42}  \rput(3.65,1.2){\scriptsize\red 124}
\cnode(4.5,2.799){2pt}{43}  \rput(4.95,2.799){\scriptsize 214}
\cnode*[linecolor=red](1.5,2.799){2pt}{46}  \rput(1.05,2.799){\scriptsize\red 234}
\cnode(2.25,4.098){2pt}{45}  \rput(2.15,4.398){\scriptsize 134}
\cnode*[linecolor=red](3.75,4.098){2pt}{44}  \rput(3.85,4.398){\scriptsize\red 314}
\ncline[linecolor=blue]{41}{42}  \ncline[linecolor=blue]{43}{42}  \ncline[linecolor=blue]{43}{44}  \ncline[linecolor=blue]{44}{45} \ncline[linecolor=blue]{45}{46} \ncline[linecolor=blue]{46}{41}
\rput(3.2,2.5){\blue $S_{3,2}^4$}

\cnode*[linecolor=red](-2.25,1.5){2pt}{11}  \rput(-2.55,1.7){\scriptsize\red 241}
\cnode(-3.75,1.5){2pt}{12}  \rput(-3.65,1.2){\scriptsize 421}
\cnode*[linecolor=red](-4.5,2.799){2pt}{13}  \rput(-4.95,2.799){\scriptsize\red 321}
\cnode(-1.5,2.799){2pt}{16}  \rput(-1.05,2.799){\scriptsize 341}
\cnode*[linecolor=red](-2.25,4.098){2pt}{15}  \rput(-2.15,4.398){\scriptsize\red 431}
\cnode(-3.75,4.098){2pt}{14}  \rput(-3.85,4.398){\scriptsize 231}
\ncline[linecolor=blue]{11}{12}  \ncline[linecolor=blue]{13}{12}  \ncline[linecolor=blue]{13}{14}  \ncline[linecolor=blue]{14}{15} \ncline[linecolor=blue]{15}{16} \ncline[linecolor=blue]{16}{11}
\rput(-3.2,2.5){\blue $S_{3,2}^1$}

\cnode(-2.25,-1.5){2pt}{31}  \rput(-2.55,-1.7){\scriptsize 413}
\cnode*[linecolor=red](-3.75,-1.5){2pt}{32}  \rput(-3.65,-1.2){\scriptsize\red 213}
\cnode(-4.5,-2.799){2pt}{33}  \rput(-4.95,-2.799){\scriptsize 123}
\cnode*[linecolor=red](-1.5,-2.799){2pt}{36}  \rput(-1.05,-2.799){\scriptsize\red 143}
\cnode(-2.25,-4.098){2pt}{35}  \rput(-2.15,-4.398){\scriptsize 243}
\cnode*[linecolor=red](-3.75,-4.098){2pt}{34}  \rput(-3.85,-4.398){\scriptsize\red 423}
\ncline[linecolor=blue]{31}{32}  \ncline[linecolor=blue]{33}{32}  \ncline[linecolor=blue]{33}{34}  \ncline[linecolor=blue]{34}{35} \ncline[linecolor=blue]{35}{36} \ncline[linecolor=blue]{36}{31}
\rput(-3.2,-2.5){\blue $S_{3,2}^3$}

\cnode*[linecolor=red](2.25,-1.5){2pt}{21}  \rput(2.55,-1.7){\scriptsize\red 132}
\cnode(3.75,-1.5){2pt}{22}  \rput(3.65,-1.2){\scriptsize 312}
\cnode*[linecolor=red](4.5,-2.799){2pt}{23}  \rput(4.95,-2.799){\scriptsize\red 412}
\cnode(1.5,-2.799){2pt}{26}  \rput(1.05,-2.799){\scriptsize 432}
\cnode*[linecolor=red](2.25,-4.098){2pt}{25}  \rput(2.15,-4.398){\scriptsize\red 342}
\cnode(3.75,-4.098){2pt}{24}  \rput(3.85,-4.398){\scriptsize 142}
\ncline[linecolor=blue]{21}{22}  \ncline[linecolor=blue]{23}{22}  \ncline[linecolor=blue]{23}{24}  \ncline[linecolor=blue]{24}{25} \ncline[linecolor=blue]{25}{26} \ncline[linecolor=blue]{26}{21}
\rput(3.2,-2.5){\blue $S_{3,2}^2$}

\ncline{11}{24}  \ncarc[arcangle=20]{42}{12} \ncline{13}{33}  \ncline{14}{21}  \ncline{15}{45} \ncline{16}{36}
\ncline{41}{34}  \ncline{43}{23}  \ncline{44}{31} \ncline{46}{26} \ncarc[arcangle=20]{32}{22} \ncline{35}{25}

\rput(0,-4.8){$S_{4,3}$}

\end{pspicture}
\caption{\label{f2}\footnotesize The $(3,2)$-star graph $S_{3,2}$ and the $(4,3)$-star graph $S_{4,3}$.}
\end{center}
\end{figure}

\begin{exap}[See Figure~\ref{f2}]
By Proposition~\ref{pps3.3}, $I_1^2=\{21,32,13\}$ and $I_2^2=\{31,12,23\}$ are
two maximum independent sets of $S_{3,2}$. The constructed two maximum independent sets of $S_{4,3}$ are:\\
$$\begin{array}{rl} I_1^3&=\{124,234,314\}\cup\{241,321,431\}\cup\{412,342,132\}\cup\{213,423,143\}\\
&=I_1^3(4)\cup I_1^3(1)\cup I_1^3(2)\cup I_1^3(3),
\end{array}$$
$$\begin{array}{rl} I_2^3&=\{134,214,324\}\cup\{341,421,231\}\cup\{312,142,432\}\cup\{413,123,243\}\\
&=I_2^3(4)\cup I_2^3(1)\cup I_2^3(2)\cup I_2^3(3).
\end{array}$$
\end{exap}

\end{document}